\documentclass[
final,nomarks
]{dmtcs-episciences}

\usepackage[utf8]{inputenc}
\usepackage{subfigure}

\usepackage{amsmath}
\usepackage{booktabs}
\usepackage{array}
\newtheorem{thm}{Theorem}[section]

\newtheorem{conjecture}[thm]{Conjecture}
\newtheorem{claim}{Claim}

\newtheorem{obs}{Obervation}[section]
\usepackage[round]{natbib}

\author{Peixue Zhao\affiliationmark{1}
  \and Fei Huang\affiliationmark{1}\thanks{Corresponding author: Fei Huang. Email: hf@zzu.edu.cn}}
\title{Rainbow vertex pair-pancyclicity of strongly edge-colored graphs}
\affiliation{
  School of Mathematics and Statistics, Zhengzhou University, Zhengzhou, Henan, People's Republic of China}
\begin{document}
\publicationdata{vol. 25:1}{2023}{13}{10.46298/dmtcs.10142}{2022-10-13; 2022-10-13; 2023-03-27}{2023-04-01}

\maketitle
\begin{abstract}
  An edge-colored graph is \emph{rainbow }if no two edges of the graph have the same color. An edge-colored graph $G^c$ is called \emph{properly colored} if every two adjacent edges of $G^c$ receive distinct colors in $G^c$. A  \emph{strongly  edge-colored}  graph is a proper edge-colored graph such that every path of length $3$ is rainbow.  We call an edge-colored graph $G^c$ \emph{rainbow vertex pair-pancyclic} if any two vertices in $G^c$ are contained in a rainbow cycle of length $\ell$ for each $\ell$ with $3 \leq \ell \leq n$. In this paper, we show that every strongly edge-colored graph $G^c$ of order $n$ with minimum degree $\delta \geq \frac{2n}{3}+1$ is rainbow vertex pair-pancyclicity.
\end{abstract}

\section{Introduction}

In this paper, we only consider finite, undirected and simple graphs. Let $G$ be a graph consisting of a vertex set $V(G)$ and an edge set $E=E(G)$. We use $d(v)$ to denote the number of edges incident with vertex $v$ in $G$ and $\delta(G)=min\{ d(v):v \in E \}$. An \emph{edge-coloring} of $G$ is a mapping $c:\ E(G) \rightarrow S$, where $S$ is a set of colors. A graph $G$ with an edge-coloring $c$ is called an  \emph{edge-colored} graph, and denoted by $G^c$. For any $e \in E(G)$, $e$ has color $k$ if $c(e)=k$. For any subset $E_1 \subseteq E$, $c(E_1)$ is the set $\{c(e):e \in E_1\}$. We use $d^c_G(v)$ (or briefly $d^c(v)$) to denote the number of different colors on the edges incident with  vertex $v$ in  $G^c$ and $\delta ^c(G)=min\{ d^c(v):v \in V(G^c) \}$. An edge-colored graph $G^c$ is called \emph{properly colored} if every two adjacent edges of $G^c$ receive distinct colors in $G^c$. Edge-colored graph $G^c$ is $rainbow$ if no two  edges of $G^c$ have the same color.  A  \emph{strongly  edge-colored graph} is a proper edge-colored graph such that every path of length $3$ is rainbow. It is clearly that $d(v)=d^c(v)$ for all $v \in V(G^c)$ in a strongly edge-colored graph $G^c$, or equivalently, for every vertex $v$ in strongly edge-colored graph $G^c$, the colors on the edges incident with $v$ are pairwise distinct.  An edge-colored graph $G^c$ is called \emph{rainbow  Hamiltonian} if $G^c$ contains a rainbow Hamiltonian cycle and \emph{rainbow  vertex(edge)-pancyclic} if every vertex (edge) in $G^c$ is contained in a rainbow cycle of length $l$ for each $l$ with $3 \leq l \leq n$. We call an edge-colored graph $G^c$ \emph{rainbow vertex  pair-pancyclic} if any two vertices in $G^c$ are contained a rainbow cycle of length $l$ for each $l$ with $3 \leq l \leq n$. further, we call a cycle $C$ \emph{l-cycle} if the length of the cycle $C$ is $l$.
For notation and terminology not defined here, we refer the reader to \cite{1}.

The classical Dirac's theorem states that every graph $G$ is Hamiltonian if $\delta(G) \geq \frac{n}{2}$. Inspired by this famous theorem, \cite{2} show that every graph $G$ of order $n$ with minimum degree $\delta \geq \frac{n+1}{2}$ is vertex-pancyclic. During the past few decades, the existence of cycles in graphs have been extensively studied in the literatures. We recommend  \cite{9,7,8,6,CHY,3,4,11,GHY,5,LHY} for more results.

For edge-colored graphs, \cite{10} proved the following asymptotic theorem about properly colored cycles.

\begin{thm}[\cite{10}] \label{thm1}
For any $\varepsilon >0$, there exists an integer $n_0$ such that every edge-colored graph $G^c$ with $n$ vertices and $\delta ^c (G) \geq (\frac{2}{3}+\varepsilon)n$ and $n \geq n_0$ contains a properly edge-colored cycle of length $l$ for all $3 \leq l \leq n$, where $\delta ^c (G)$ is the minimum number of distinct colors of edges incident with a vertex in $G^c$.
\end{thm}

\cite{12} considered the existence of rainbow Hamiltonian cycles in strongly edge-colored graph and proposed the following two conjectures.

\begin{conjecture}[\cite{12}]\label{conject1}
Every strongly edge-colored graph $G^c$ with $n$ vertices and degree at least $\frac{n+1}{2}$ has a rainbow Hamiltonian cycle.
\end{conjecture}

\begin{conjecture}[\cite{12}]\label{conject2}
Every strongly edge-colored graph $G^c$ with $n$ vertices and degree at least $\frac{n}{2}$ has a rainbow Hamiltonian path.
\end{conjecture}

To support the above two conjectures, they presented the following theorem.

\begin{thm}[\cite{12}] \label{thm4}
Let $G^c$ be a strongly edge-colored graph with minimum degree $\delta$, if $\delta \geq \frac{2 |G|}{3}$, then $G^c$ has a rainbow Hamiltonian cycle.
\end{thm}

\cite{13} showed that every strongly edge-colored graph $G^c$  on $n$ vertices is rainbow vertex-pancyclic if $\delta \geq \frac{2n}{3}$. \cite{14} further considered the rainbow edge-pancyclicity of  strongly edge-colored graphs and proposed the following theorem.

\begin{thm}[\cite{14}]\label{ll} Let $G^c$ be a strongly edge-colored graph on  $n$ vertices. If $\delta(G^c) \geq \frac{2n+1}{3}$, then $G^c$ is rainbow edge-pancyclic. Furthermore, for every edge $e$ of $G^c$, one can find a rainbow $l$-cycle containing $e$ for each $l$ with $3 \leq l \leq n$ in polynomial time.
\end{thm}

In this paper, we consider the rainbow vertex pair-pancyclicity of strongly edge-colored graph. Our main result is as follows.

\begin{thm} \label{thmme}                                                                                    
Let $G^c$ be a strongly edge-colored graph with $n$ vertices and minimum degree $\delta$. If $\delta \geq \frac{2n}{3}+1$, then $G^c$ is rainbow vertex pair-pancyclicity.
\end{thm}

\section{Proof of Theorem \ref{thmme}}                                                   

First, we introduce some useful notations. Given a rainbow cycle $C$ in graph $G^c$, a color $s$ is called a \emph{$C$-color} (resp., \emph{$\widetilde{C}$-color}) if $s \in c(C)$ (resp., $s \notin c(C)$). Correspondingly, we call an edge $e$ a \emph{$C$-color edge} (resp., \emph{$\widetilde{C}$-color edge}) if $c(e) \in c(C)$ (resp., $c(e) \notin c(C)$). Two adjacent vertices $u$ and $v$ are called \emph{$C$-adjacent} (resp., \emph{$\widetilde{C}$-adjacent}) if $c(uv) \in c(C)$ (resp., $c(uv) \notin c(C)$). For two disjoint adjacent subsets $V_1$ and $V_2$ of $V(G)$, let $E(V_1,V_2)$ denote the set of edges between $V_1$ and $V_2$. We denote the subsets of $E(V_1,V_2)$ consisting of the $C$-color edges (resp., $\widetilde{C}$-color edges) by $E_C(V_1,V_2)$ (resp., $E_{\widetilde{C}}(V_1,V_2)$). Similarly, for two subgraphs $H_1$ and $H_2$, we denote the set of $C$-color edges (resp., $\widetilde{C}$-color edges) between $V(H_1)$ and $V(H_2)$ by $E_C(H_1,H_2)$ (resp., $E_{\widetilde{C}}(H_1,H_2)$). For any two vertices $v_i$ and $v_j$ of cycle $C=v_1v_2\ldots v_lv_1$,  we identify the two subscripts $i$ and $j$ if $i\equiv j\ (mod \ l)$. Let $v_iC^+v_j$ be the path $v_iv_{i+1}...v_{j-1}v_j$ and $v_iC^-v_j$ the path $v_iv_{i-1}...v_{j+1}v_j$, respectively. For any vertex $v \in V(G^c)$, let $CN(v)$ be the set of colors used by the edges incident with $v$.

From the definition of strongly edge-coloring, we can easily get the following observation.

\begin{obs}\label{ob1}
Each cycle of length at most 5 in a strongly edge-colored graph is rainbow.
\end{obs}

\noindent\textbf{Proof of Theorem \ref{thmme}:} Recall that the colors on the edges incident with $v$ are pairwise distinct for each vertex $v$ of a strongly edge-colored graph. So we do not distinguish the colors of adjacent edges in the following.
 If $n \leq 8$, $G$ is complete since $\delta \geq \frac{2n}{3}+1$, and so the result clearly holds. Thus we suppose that $n \geq 9$. Let $a$ and $b$ be two arbitrary vertices of $G$. If $a$ and $b$ are adjacent, then $a$ and $b$ are contained in a rainbow cycle of length $l$ for each $l$ with $3 \leq l \leq n$ from Theorem \ref{ll}. So we  consider that $a$ and $b$ are not adjacent. Since  $\delta \geq \frac{2n}{3}+1$, we have that $a$ and $b$ are contained in a 4-cycle which is rainbow from Observation \ref{ob1}.  Suppose to the contrary that the result is not true. Then there is an integer $l$ with $4\le l\le n-1$ such that there is a rainbow $l$-cycle containing $a$ and $b$, but there is no rainbow $(l+1)$-cycle containing  both $a$ and $b$.  Let $C:=v_1v_2\ldots v_lv_1$ be a rainbow $l$-cycle containing $a$ and $b$.

Without loss of generality, we assume that $c(v_iv_{i+1})=i$ for $1 \leq i \leq l$.  For $1\le i\le l$, let $N_i$ be the set of the vertices of $C$ which are adjacent to $v_i$, that is, $N_i=N(v_i)\cap V(C)$. We then proof the following claim.
\begin{claim} \label{claim1}
 $l \geq \frac{n+12}{3}$. In particular, $l \geq 7$ when $n\geq 9$.
\end{claim}

\noindent\textbf{Proof.}   Since $G^c$ is strongly edge-colored, for any $v_j \in N_1$, the color $j$ does not occur in $CN(v_1)$. So the number of $C$-colors  not contained in $CN(v_1)$ is at least $|N_1|-1$, and therefore, the number of $C$-colors contained in $CN(v_1)$ is at most $l-(|N_1|-1)$. Since $1$ and $l$ are $C$-colors in $CN(v_1)$, we have that the number of $C$-colors contained in $E(v_1,V(G) \setminus V(C))$ is at most $l-(|N_1|-1)-2=l-|N_1|-1$. Hence, we have  $|E_C(v_1,V(G) \setminus V(C))| \leq l-|N_1|-1$. Since $|E(v_1,V(G) \setminus V(C))| \geq \delta - |N_1|$,  we have  that
\begin{equation*}
\begin{split}
|E_{\widetilde{C}}(v_1,V(G) \setminus V(C))| &= |E(v_1,V(G) \setminus V(C))|-|E_C(v_1,V(G) \setminus V(C))|\\
                                             &\geq (\delta -|N_1|)-(l-|N_1|-1)\\
                                             &= \delta -l+1.
\end{split}
\end{equation*}
Similarly, we can also deduce that $|E_{\widetilde{C}}(v_i,V(G) \setminus V(C))| \geq \delta -l+1$ for all $1 \leq i\leq l$. For any two vertices $v_i$ and $v_{i+1}$ with $1 \leq i \leq l$, if there exists a vertex $w \in V(G) \setminus V(C)$ such that both $v_iw$ and $v_{i+1}w$ are $\widetilde{C}$-color edges, then both  $a$ and $b$ are contained in a rainbow $(l+1)$-cycle $C':=v_iwv_{i+1}C^+v_i$, a contradiction. Thus, for any common neighbor $w \in V(G) \setminus V(C)$ of $v_i$ and $v_{i+1}$, either $v_iw$ or $v_{i+1}w$ is not a $\widetilde{C}$-color edge. Then we have that $|E_{\widetilde{C}}(v_i,w)|+|E_{\widetilde{C}}(v_{i+1},w)| \leq 1$. Therefore, we have
\[
n \geq |E_{\widetilde{C}}(v_i,V(G)\setminus V(C))|+|E_{\widetilde{C}}(v_{i+1},V(G)\setminus V(C))|+l \geq 2(\delta - l + 1)+l=2\delta -l+2.
\]
Hence,
\begin{equation*}
l\geq 2\delta -n +2 \geq 2 \cdot (\frac{2n}{3}+1)-n+2=\frac{n+12}{3}.
\label{eq1}
\end{equation*}                                                                                     
This completes the claim.
$\hfill\square$\\

Let $H=K_k$ be the maximal rainbow complete graph in $G^c[V(G) \setminus V(C)]$ such that every edge in $H$ is $\widetilde{C}$-colored, and let $R=G^c[V(G)-(V(C) \cup V(H))]$. It is clearly that for any $w \in V(H)$, if there is a vertex $v_i \in V(C)$ such that $v_iw$ is a $\widetilde{C}$-color edge, then $c(v_iw) \notin c(H)$ since $G^c$ is a strongly edge-colored graph.

For two $\widetilde{C}$-color edges $v_iw_1$ and $v_jw_2$ with $w_1,w_2 \in V(H)$ and $1 \leq i < j \leq l$, if $w_1=w_2$ and $j-i=1$, we say  $v_iw_1$ and $v_jw_2$ are  \emph{forbidden pair  of type  1};  if $w_1 \neq w_2$,  both $a$ and $b$ are contained in $v_iC^-v_j$, and $2 \leq j-i \leq k$, we say $v_iw_1$ and $v_jw_2$ are  \emph{forbidden pair  of type  2}. Clearly, if $E_{\widetilde{C}}(C,H)$ has a forbidden pair of type 1, then there exists a rainbow $(l+1)$-cycle $C':=v_iw_1v_jC^+v_i$ containing both $a$ and $b$, and if $E_{\widetilde{C}}(C,H)$ has a forbidden pair of type 2, then there exist a rainbow $(l+1)$-cycle $C':=v_iw_1Hw_2v_jC^+v_i$ containing both $a$ and $b$, where $w_1Hw_2$ is a path of length $|E(v_iC^+v_j)|-1$ with endpoints $w_1$ and $w_2$ in $H$.
\begin{claim} \label{claim2}
$k \geq 3.$
\end{claim}                                                       

\noindent\textbf{Proof.} For each $w \in V(H)$, let
\begin{gather*}
\widetilde{s}_w=|E_{\widetilde{C}}(w,C)|, s_w=|E_C(w,C)|, \\
\widetilde{t}_w=|E_{\widetilde{C}}(w,R)|, t_w=|E_C(w,R)|.
\end{gather*}
We have
\begin{equation}
\widetilde{s}_w + s_w + \widetilde{t}_w +t_w +(k-1) \geq \delta.
\label{eq2}
\end{equation}                                                                                      
If there is an integer $i$ with $1 \leq i \leq l$ such that $v_iw \in E(G^c)$, then the colors $i-1$ and $i$ can not appear in $CN(w)$. Thus the number of $C$-colors not contained in $CN(w)$ is at least $\widetilde{s}_w + s_w$, which implies that
\[
s_w + t_w \leq l-(\widetilde{s}_w + s_w),
\]
and so, we have
\begin{equation}
\widetilde{s}_w +2s_w +t_w \leq l.
\label{eq3}
\end{equation}                                                                                      
Let $v_{i_1},v_{i_2},...,v_{i_{\widetilde{s}_w}}$ be the vertices on $C$ which are $\widetilde{C}$-adjacent to $w$. Without loss of generality, we suppose that $1 \leq i_1 < i_2 < ... < i_{\widetilde{s}_w} \leq l$. Then $i_{j+1} - i_j \geq 2$ for each $1 \leq j \leq \widetilde{s}_w-1$ and $i_{\widetilde{s}_w}-i_1 \leq l-2$.  Let $I = \{ i_1-1,i_1,i_2-1,i_2,...,i_{\widetilde{s}_w}-1,i_{\widetilde{s}_w} \}$. Clearly, we have $|I|=2 \widetilde{s}_w$ and $I \cap CN(w) = \phi$. Thus, we can deduce that
\begin{equation}
2\widetilde{s}_w + s_w + t_w = |I| + s_w + t_w \leq l.
\label{eq4}
\end{equation}                                                                                      
Since $|V(R)|=n-l-k$, we have $
t_w + \widetilde{t}_w \leq n-l-k.
$
Together with inequalities (\ref{eq3}) and (\ref{eq4}), we have
\[
3\widetilde{s}_w + 3s_w + 3t_w + \widetilde{t}_w  \leq l+l+n-l-k =n+l-k.
\]
Let
\[
\widetilde{S}= \sum_{w \in V(H)}\widetilde{s}_w,
S= \sum_{w \in V(H)}s_w,
\widetilde{T}= \sum_{w \in V(H)}\widetilde{t}_w,
T= \sum_{w \in V(H)}t_w.
\]
Then,
\begin{equation}
3\widetilde{S} + 3S + 3T + \widetilde{T} \leq k(n+l-k).
\label{eq5}
\end{equation}                                                                                      
Since $k$ is maximal, each vertex of $R$ has at most $k-1$ number of $\widetilde{C}$-color edges to $H$, which implies that
\begin{equation}
\widetilde{T}= \sum_{w \in V(H)}\widetilde{t}_w \leq (k-1)(n-l-k).
\label{eq6}
\end{equation}                                                                                      
Recall that $w \in V(H)$. By (\ref{eq2}) and the arbitrariness of $w$, we have
\begin{equation}
\begin{split}
k \delta &\leq \sum_{w \in V(H)}(\widetilde{s}_w + s_w + \widetilde{t}_w +t_w +(k-1)) \\
         &=\widetilde{S} + S + \widetilde{T} + T +k(k-1).\label{eq7}
\end{split}
\end{equation}                                                                                      
Combining inequalities (\ref{eq5}), (\ref{eq6}) and (\ref{eq7}), we can get the following inequality
\begin{equation*}
\begin{split}
3k \delta &\leq 3\widetilde{S} + 3S + 3T + 3\widetilde{T} + 3k(k-1) \\
          &\leq k(n+l-k) + 2(k-1)(n-l-k) + 3k(k-1) \\
          &\leq n(3k-2) + l(2-k) - k.
\end{split}
\end{equation*}
If $k=1$, then $l > n$, a contradiction. If $k=2$, then $\delta \leq \frac{2n-1}{3}$, again a contradiction. So we have $k \geq 3$. Claim 2 follows.
$\hfill\square$\\

Since $H$ is a rainbow complete graph, we can deduce that
\begin{equation}
S+T \leq l.
\label{eq8}
\end{equation}
\begin{claim} \label{claim3}
 $\widetilde{S} \geq l+1$.
 \end{claim}

\noindent\textbf{Proof.} Suppose, by way of contradiction, that $\widetilde{S} \leq l$. Combining with inequality (\ref{eq7}), we can get that  $$k \delta \leq \widetilde{S} + S + \widetilde{T} + T +k(k-1)
         \leq l + l + (k-1)(n-l-k) + k(k-1),$$
which implies that $k(n-l- \delta) \geq n-3l$. Since $\delta \ge \frac{2n}{3}+1$ and $l\ge \frac{n+12}{3}$ from Claim \ref{claim1}, we  have  $n-l- \delta \leq 0$. Thus we have $3(n-l- \delta) \geq k(n-l- \delta) \geq n-3l$ from Claim \ref{claim2},  and therefore
$\delta \leq \frac{2n}{3}$, a contradiction. Claim 3 follows.
$\hfill\square$\\

Without loss of generality, we suppose that $a=v_1$ and $b=v_m$, where $2 \leq m \leq l-1$, and let  $P^1=aC^+b$. Then we design an algorithm to generate a sequence of disjoint sub-paths $P^1_1,P^1_2,...,P^1_{h_1}$ of $C$ respect to $P^1$ and $H$.\\

{\noindent
\begin{tabular}{m{13cm}}

\textbf{Algorithm AI} \\
\toprule
\textbf{Input:} a strongly edge-colored graph $G^c$, a rainbow cycle $C=v_1v_2\ldots v_lv_1$, a path $P^1 = v_1v_2...v_m$ and a rainbow complete subgraph $H=K_k$ of $G^c-V(C)$. \\
\textbf{Output:} a sequence of  disjoint paths $P^1_1,P^1_2,...,P^1_{h_1}$ such that $P^1_i$ is a subgraph of $C$. \\
1: \textbf{Set} $i=1$ \\
2: \textbf{While} $V(P^1) \neq \phi$ \textbf{do} \\
\quad \textbf{If} $E_{\widetilde{C}}(P^1,H) = \phi$ \\
\quad \quad \textbf{stop} \\
\quad \textbf{Else} \textbf{Set} $d$ be the smallest subscript such that $E_{\widetilde{C}}(v_d,H) \neq \phi$ \\
\quad \quad \textbf{If} $d+k\geq m$ \textbf{then} \\
\quad \quad \quad \textbf{Set} $P^1_i=v_dv_{d+1}...v_m$ \\
\quad \quad \quad \textbf{stop} \\
\quad \quad \textbf{Else}\quad \textbf{If} $|E_{\widetilde{C}}(v_d,H)| \geq 2$ \textbf{then} \\
\quad \quad \quad \quad \quad \quad \textbf{Set} $P^1_i=v_dv_{d+1}...v_{d+k}$ \\
\quad \quad \quad \quad \quad \textbf{If} $|E_{\widetilde{C}}(v_d,H)| =1$ \textbf{then} \\
\quad \quad \quad \quad \quad \quad \textbf{Set} $P^1_i=v_dv_{d+1}...v_{d+k+1}$ \\
\quad \quad \textbf{Set} $P^1=P^1 \setminus P^1_i$ \\
\quad \quad \textbf{Set} $i=i+1$ \\
3: \textbf{return} $P^1_1,P^1_2,...,P^1_{h_1}$ \\
\bottomrule
\end{tabular}
}                                                                                          

\begin{claim} \label{claim4}
 $|E_{\widetilde{C}}(P^1_i,H)| \leq |V(P^1_i)|-1$ for any  $1 \leq i \leq h_1-1$,  $|E_{\widetilde{C}}(P^1_{h_1},H)| \leq k$ if $|V(P^1_{h_1})| \in \{ 1,2 \}$, and $|E_{\widetilde{C}}(P^1_{h_1},H)| \leq k + 1$ if $3 \leq |V(P^1_{h_1})| \leq k+1$.
\end{claim}

\noindent\textbf{Proof.}
For $1 \leq i \leq h_1-1$, we distinguish the following two cases.\\
\textbf{Case 1.} $|E_{\widetilde{C}}(v_d,H)| \geq 2.$ Then we have $P^1_i=v_dv_{d+1}...v_{d+k}$. Let $w_1$ and $w_2$ be two vertices in $H$ such that $v_dw_1,v_dw_2 \in E_{\widetilde{C}}(v_d,H)$. Since there exist no forbidden pairs of type 1 for any vertex $w \in V(H)$, then we have $|E_{\widetilde{C}}(v_d,H)| + |E_{\widetilde{C}}(v_{d+1},H)| \leq k$.  For any  $j$ with $d+2 \leq j \leq d+k$, if  $w_1$ and $v_j$ are $\widetilde{C}$-adjacent, then $v_jw_1$  and $v_dw_2$ form a forbidden pair of type 2; if  $w_2$ and $v_j$ are
$\widetilde{C}$-adjacent, then  $v_jw_2$  and $v_dw_1$ form a forbidden pair of type 2; if $v_j$ and $w$ are $\widetilde{C}$-adjacent for some $w$ with $w \neq w_1$ and $w \neq w_2$, then $v_jw$ and $v_dw_1$  form a forbidden pair of type 2. Therefore, we have $|E_{\widetilde{C}}(v_j,H)|=0$. Thus,
\begin{equation*}
\begin{split}
|E_{\widetilde{C}}(P^1_i,H)| &= \sum^{d+k}_{j=d} |E_{\widetilde{C}}(v_j,H)| \\
                             &= |E_{\widetilde{C}}(v_d,H)| + |E_{\widetilde{C}}(v_{d+1},H)| \\
                             &\leq k \\
                             &= |V(P^1_i)|-1.
\end{split}
\end{equation*}\\
\textbf{Case 2.} $|E_{\widetilde{C}}(v_d,H)| = 1.$  Then we have $P^1_i=v_dv_{d+1}...v_{d+k+1}$. Let $w_1$ be a vertex in $H$ such that $v_dw_1 \in E_{\widetilde{C}}(v_d,H)$. We further distinguish the following three cases.\\
\textbf{Case 2.1.} $|E_{\widetilde{C}}(v_{d+1},H)|=0.$ For any $w \in V(H) \setminus \{ w_1 \}$, we have that $v_{j}$ and $w$ cannot be $\widetilde{C}$-adjacent for any  $d+2\le j\le d+k+1$ since otherwise  $v_jw$ and $v_dw_1$ form a forbidden pair of type 2. Thus, we have $|E_{\widetilde{C}}(v_j,H)| \leq 1$ and $\sum ^{d+k+1}_{j=d+2} |E_{\widetilde{C}}(v_j,H)| \leq k-1$. Therefore, 
\begin{equation*}
\begin{split}
|E_{\widetilde{C}}(P^1_i,H)| &= \sum^{d+k+1}_{j=d} |E_{\widetilde{C}}(v_j,H)|  \\
                             &=|E_{\widetilde{C}}(v_d,H)| + |E_{\widetilde{C}}(v_{d+1},H)| + \sum ^{d+k+1}_{j=d+2} |E_{\widetilde{C}}(v_j,H)|  \\
                             &\leq 1+0+(k-1) \\
                             &=k \\
                             &\le |V(P^1_i)|-1.
\end{split}
\end{equation*}\\
\textbf{Case 2.2.} $|E_{\widetilde{C}}(v_{d+1},H)|=1.$ Let $w_2$ be a vertex in $H$ such that $v_{d+1}w_2 \in E_{\widetilde{C}}(v_d,H)$. Clearly, $w_1 \neq w_2$. If $v_{d+2}$ and $w_2$ are $\widetilde{C}$-adjacent, we have that $v_{d+2}w_2$ and $v_dw_1$ form a forbidden pair of type 2, a contradiction. If $v_{d+2}$ and $w$ are $\widetilde{C}$-adjacent for some $w\in V(H)$ with $w \neq w_1$ and $w \neq w_2$, then $v_{d+2}w$ and $v_dw_1$ form a forbidden pair of type 2, again a contradiction. So, $|E_{\widetilde{C}}(v_{d+2},H)| \leq 1.$
For any  $j$ with $d+3 \leq j \leq d+k+1$, if  $w_1$ and $v_j$ are $\widetilde{C}$-adjacent, then $v_jw_1$  and $v_{d+1}w_2$ form a forbidden pair of type 2; if  $w_2$ and $v_j$ are
$\widetilde{C}$-adjacent, then  $v_jw_2$  and $v_dw_1$ form a forbidden pair of type 2; if $v_j$ and $w$ are $\widetilde{C}$-adjacent for some $w\in V(H)$ with $w \neq w_1$ and $w \neq w_2$, then $v_jw$ and $v_dw_1$  form a forbidden pair of type 2. We obtain a contradiction in the above three cases, and therefore, we have $\sum ^{d+k+1}_{j=d+3} |E_{\widetilde{C}}(v_{j},H)|=0.$ Therefore,

\begin{equation*}
\begin{split}
|E_{\widetilde{C}}(P^1_i,H)| &= \sum^{d+k+1}_{j=d} |E_{\widetilde{C}}(v_j,H)|  \\
                             &= |E_{\widetilde{C}}(v_d,H)| +  |E_{\widetilde{C}}(v_{d+1},H)| + |E_{\widetilde{C}}(v_{d+2},H)| + \sum ^{d+k+1}_{j=d+3} |E_{\widetilde{C}}(v_j,H)|   \\
                             &\leq 1+1+1+0 \\
                             &\le k \\
                             &\le |V(P^1_i)|-1.
\end{split}
\end{equation*}\\
\textbf{Case 2.3.} $|E_{\widetilde{C}}(v_{d+1},H)| \geq 2.$ Let $Q^1_i =P^1_i \setminus \{ v_d \}= v_{d+1}v_{d+2}...v_{d+k+1}$.  Similar to the discussion of  Case 1, we have that $|E_{\widetilde{C}}(Q^1_i,H)| \leq |V(Q^1_i)|-1 = (k+1)-1 = k$. Thus, $|E_{\widetilde{C}}(P^1_i,H)| = |E_{\widetilde{C}}(v_d,H)| + |E_{\widetilde{C}}(Q^1_i,H)| \leq 1+k=|V(P^1_i)|-1$.

Then we analysis the value of $|E_{\widetilde{C}}(P^1_{h_1},H)|$. If $|V(P^1_{h_1})| = 1$, the inequality $|E_{\widetilde{C}}(P^1_{h_1},H)| \leq k$ clearly holds. If $|V(P^1_{h_1})| = 2$, that is, $P^1_{h_1}=v_dv_{d+1}$,  we have $|E_{\widetilde{C}}(v_d,H)| + |E_{\widetilde{C}}(v_{d+1},H)| \leq k$ since $v_d$ and $v_{d+1}$ are adjacent. Therefore, $|E_{\widetilde{C}}(P^1_{h_1},H)|= E_{\widetilde{C}}(v_d,H)| + |E_{\widetilde{C}}(v_{d+1},H)| \leq k$. If $3 \leq |V(P^1_{h_1})| \leq k+1$, we have $|E_{\widetilde{C}}(P^1_{h_1},H)| \leq k$ when  $|E_{\widetilde{C}}(v_d,H)| \geq 2$  by the similar analysis of the above Case 1 (taking $m$ as $d+k$), and $|E_{\widetilde{C}}(P^1_{h_1},H)| \leq k+1$ when $|E_{\widetilde{C}}(v_d,H)| =1$ by the similar analysis of the above Case 2  (taking $m$ as $d+k+1$).  The proof is thus completed. $\hfill\square$\\

Let $P^2=aC^-b$.  Then we design another algorithm to generate a sequence of disjoint sub-paths  $P^2_1,P^2_2,...,P^2_{h_2}$ of $C$ respect to  $P^2$ and  $H$ in the following.\\

{\noindent
\begin{tabular}{m{13cm}}

\textbf{Algorithm AII} \\
\toprule
\textbf{Input:} a strongly edge-colored graph $G$, a rainbow cycle $C=v_1v_2\ldots v_lv_1$, $P^2=aC^-b = v_{l+1}v_lv_{l-1}...v_m$ and a rainbow complete subgraph $H=K_k$ of $G^c-V(C)$.  \\
\textbf{Output:} a sequence of  disjoint paths $P^2_1,P^2_2,...,P^2_{h_2}$ such that $P^2_i$ is a subgraph of $C$. \\
1: \textbf{Set} $i=1$ \\
2: \textbf{While} $V(P^2) \neq \phi$ \textbf{do} \\
\quad \textbf{If} $E_{\widetilde{C}}(P^2,H) = \phi$ \\
\quad \quad \textbf{stop} \\
\quad \textbf{Else} \textbf{Set} $d$ be the biggest subscript for which $E_{\widetilde{C}}(v_d,H) \neq \phi$ \\
\quad \quad \textbf{If} $d-k \leq m$ \textbf{then} \\
\quad \quad \quad \textbf{Set} $P^2_i=v_dv_{d-1}...v_m$ \\
\quad \quad \quad \textbf{stop} \\
\quad \quad \textbf{Else}\quad \textbf{If} $|E_{\widetilde{C}}(v_d,H)| \geq 2$ \textbf{then} \\
\quad \quad \quad \quad \quad \quad \textbf{Set} $P^2_i=v_dv_{d-1}...v_{d-k}$ \\
\quad \quad \quad \quad \quad \textbf{If} $|E_{\widetilde{C}}(v_d,H)| =1$ \textbf{then} \\
\quad \quad \quad \quad \quad \quad \textbf{Set} $P^2_i=v_dv_{d-1}...v_{d-k-1}$ \\
\quad \quad \textbf{Set} $P^2=P^2 \setminus P^2_i$ \\
\quad \quad \textbf{Set} $i=i+1$ \\
3: \textbf{return} $P^2_1,P^2_2,...,P^2_{h_2}$ \\
\bottomrule
\end{tabular}
}                                                                                        
\\[4mm]

Similar to Claim \ref{claim4}, we can get the following Claim.
\begin{claim} \label{claim5}
$|E_{\widetilde{C}}(P^2_i,H)| \le |V(P^2_i)|-1$ for all $1 \leq i \leq h_2-1$, $|E_{\widetilde{C}}(P^2_{h_2},H)| \leq k$ if $|V(P^2_{h_2})| \in \{ 1,2 \}$ and $|E_{\widetilde{C}}(P^2_{h_2},H)| \leq k + 1$ if $3 \leq |V(P^2_{h_2})| \leq k+1$.
\end{claim}

According to the above claims, we have
\begin{equation}
\begin{split}
|E_{\widetilde{C}}(C,H)| &=  |E_{\widetilde{C}}(aC^+b,H)| + |E_{\widetilde{C}}(aC^-b,H)| - |E_{\widetilde{C}}(a,H)| - |E_{\widetilde{C}}(b,H)| \\
                         &\leq \sum^{h_1-1}_{i=1} |V(P^1_i)| - (h_1-1) + |E_{\widetilde{C}}(P^1_{h_1},H)| \\
                         &{\quad}+ \sum^{h_2-1}_{i=1} |V(P^2_i)| - (h_2-1) + |E_{\widetilde{C}}(P^2_{h_2},H)|  \\
                         &{\quad}- |E_{\widetilde{C}}(a,H)| - |E_{\widetilde{C}}(b,H)| \\
                         &\leq [ l-|V(P^1_{h_1})|-|V(P^2_{h_2})|+1 ] - (h_1+h_2) +2 \\
                         &{\quad}+ |E_{\widetilde{C}}(P^1_{h_1},H)| + |E_{\widetilde{C}}(P^2_{h_2},H)| - |E_{\widetilde{C}}(a,H)| - |E_{\widetilde{C}}(b,H)| \\
                         &= l- ( |V(P^1_{h_1})|+|V(P^2_{h_2})| ) - (h_1+h_2) +3 \\
                         &{\quad}+ |E_{\widetilde{C}}(P^1_{h_1},H)| + |E_{\widetilde{C}}(P^2_{h_2},H)| - |E_{\widetilde{C}}(a,H)| - |E_{\widetilde{C}}(b,H)|.
\label{eq9}
\end{split}
\end{equation}                                                                                      

\begin{claim} \label{claim6}                                                                                
 $\widetilde{S} \leq l+2k-4$.
\end{claim}

\noindent\textbf{Proof.}                                                                                      
We show that $\widetilde{S} \leq max \{ 2k+2,l+k-1,l+2k-4 \}$, which implies $\widetilde{S} \leq l+2k-4$ since $l \geq 7$ from Claim \ref{claim1} and $k \geq 3$ from Claim \ref{claim2}.

Let $h=h_1+h_2$. By symmetry, we suppose $h_1 \ge h_2$ and $|V(P^1_{h_1})| \ge |V(P^2_{h_2})|$. From Claim \ref{claim3}, we have $h \geq 1$. Then we proceed our proof by distinguishing the following four cases. \\
\textbf{Case 1.} $h_1=1$ and $h_2=0$. From Algorithm AII, we have $E_{\widetilde{C}}(aC^-b,H) = \phi$. Thus, $E_{\widetilde{C}}(a,H)= \phi$ and $E_{\widetilde{C}}(b,H)= \phi$. From Algorithm AI, we have $|V(P^1_{h_1})| \geq 2$. If $|V(P^1_{h_1})| = 2$, let $u$ be the vertex distinct from $b$ in $C$  such that $E_{\widetilde{C}}(u,H) \neq \phi$. Thus we have  $\widetilde{S} = |E_{\widetilde{C}}(u,H)| \leq k < 2k+2$.  If $|V(P^1_{h_1})| \geq 3$, from Claim \ref{claim4}, we have  $\widetilde{S} = E_{\widetilde{C}}(P^1_{h_1},H) \leq k+1 < 2k+2$. The claim follows.\\
\textbf{Case 2.} $h_1 \geq 2$ and $h_2=0$. From Algorithm AI and AII, we have $E_{\widetilde{C}}(a,H)= \phi$,  $E_{\widetilde{C}}(b,H)= \phi$ and $|V(P^1_{h_1})| \geq 2$. If $|V(P^1_{h_1})| = 2$, since $E_{\widetilde{C}}(b,H)= \phi$, we have
$
|E_{\widetilde{C}}(P^1_{h_1},H)| + |E_{\widetilde{C}}(P^2_{h_2},H)| - |E_{\widetilde{C}}(b,H)|=|E_{\widetilde{C}}(P^1_{h_1},H)| \leq k.
$
Applying inequality (\ref{eq9}), we have $\widetilde{S} \leq l-2-2+3+k+0=l+k-1$. If $|V(P^1_{h_1})| \geq 3$, from Claim \ref{claim4}, we have
$
|E_{\widetilde{C}}(P^1_{h_1},H)| + |E_{\widetilde{C}}(P^2_{h_2},H)| - |E_{\widetilde{C}}(b,H)|=|E_{\widetilde{C}}(P^1_{h_1},H)| \leq k+1.
$
Thus, by inequality (\ref{eq9}), we have $\widetilde{S} \leq l-3-2+3+k+1+0=l+k-1$. The claim follows.\\
\textbf{Case 3.}  $h_1=1$ and $h_2=1$.
By Claim \ref{claim4} and \ref{claim5}, if $|V(P^1_{h_1})| \in \{ 1,2 \}$ and $|V(P^2_{h_2})| \in \{ 1,2 \}$, we have $\widetilde{S} \leq |E_{\widetilde{C}}(P^1_{h_1},H)| + |E_{\widetilde{C}}(P^2_{h_2},H)| \leq 2k <2k+2$. If $|V(P^1_{h_1})| \geq 3$ and $|V(P^2_{h_2})| \in \{ 1,2 \}$, we have $\widetilde{S} \leq |E_{\widetilde{C}}(P^1_{h_1},H)| + |E_{\widetilde{C}}(P^2_{h_2},H)| \leq 2k+1 <2k+2$. If $|V(P^1_{h_1})| \geq 3$ and $|V(P^2_{h_2})| \geq 3$, we have $\widetilde{S} \leq |E_{\widetilde{C}}(P^1_{h_1},H)| + |E_{\widetilde{C}}(P^2_{h_2},H)| \leq 2k+2$. The claim holds.\\
\textbf{Case 4.}  $h \geq 3$ and $h_2\ge 1$. We consider the following six cases.\\
\textbf{Case 4.1.} $|V(P^1_{h_1})| = 1$ and $|V(P^2_{h_2})| = 1$. It is clearly that
\[      V(P^1_{h_1})  = V(P^2_{h_2}) = \{ b \}       \]
and
\[         |E_{\widetilde{C}}(P^1_{h_1},H)| + |E_{\widetilde{C}}(P^2_{h_2},H)| - |E_{\widetilde{C}}(b,H)| = |E_{\widetilde{C}}(b,H)| \leq k.          \]
By inequality (\ref{eq9}), we have
\begin{equation*}
\widetilde{S} = |E_{\widetilde{C}}(C,H)| \leq l-2-3+3+k+0=l+k-2<l+k-1.
\end{equation*}
\textbf{Case 4.2.} $|V(P^1_{h_1})| = 2$ and $|V(P^2_{h_2})| = 1$.
It is clearly that $V(P^2_{h_2})  =  \{ b \}$.
From Claim \ref{claim4}, we have
\[         |E_{\widetilde{C}}(P^1_{h_1},H)| + |E_{\widetilde{C}}(P^2_{h_2},H)| - |E_{\widetilde{C}}(b,H)| = |E_{\widetilde{C}}(P^1_{h_1},H)| \leq k.          \]
By inequality (\ref{eq9}) and $h \geq 3$, we have
\begin{equation*}
\widetilde{S} \leq l-3-3+3+k+0 =l+k-3 < l+k-1.
\end{equation*}
\textbf{Case 4.3.} $|V(P^1_{h_1})| \geq 3$ and $|V(P^2_{h_2})| =1$. It is clearly that
$     V(P^2_{h_2})  =  \{ b \}$.
From Claim \ref{claim4}, we have
\[         |E_{\widetilde{C}}(P^1_{h_1},H)| + |E_{\widetilde{C}}(P^2_{h_2},H)| - |E_{\widetilde{C}}(b,H)| = |E_{\widetilde{C}}(P^1_{h_1},H)| \leq k+1.          \]
By inequality (\ref{eq9}) and $h \geq 3$, we have
\begin{equation*}
\widetilde{S} = |E_{\widetilde{C}}(C,H)| \leq l-4-3+3+k+1+0 =l+k-3<l+k-1.
\end{equation*}
\textbf{Case 4.4.} $|V(P^1_{h_1})| = 2$ and $|V(P^2_{h_2})| = 2$.
From  Claim \ref{claim4} and \ref{claim5}, we have
\[         |E_{\widetilde{C}}(P^1_{h_1},H)| + |E_{\widetilde{C}}(P^2_{h_2},H)| - |E_{\widetilde{C}}(b,H)| \leq 2k.          \]
By inequality (\ref{eq9}) and $h \geq 3$, we have
\begin{equation*}
\widetilde{S} = |E_{\widetilde{C}}(C,H)| \leq l-4-3+3+2k+0=l+2k-4<l+k-1.
\end{equation*}
\textbf{Case 4.5.} $|V(P^1_{h_1})| \geq 3$ and $|V(P^2_{h_2})| =2$. It is clearly that
\[         |E_{\widetilde{C}}(P^1_{h_1},H)| + |E_{\widetilde{C}}(P^2_{h_2},H)| - |E_{\widetilde{C}}(b,H)| \leq k+k+1 = 2k+1.          \]
By inequality (\ref{eq9}) and $h \geq 3$, we have
\begin{equation*}
\widetilde{S} = |E_{\widetilde{C}}(C,H)| \leq l-5-3+3+2k+1+0=l+2k-4.
\end{equation*}
\textbf{Case 4.6.} $|V(P^1_{h_1})| \geq 3$ and $|V(P^2_{h_2})| \geq 3$.
From Claim \ref{claim4} and \ref{claim5}, we have
\[         |E_{\widetilde{C}}(P^1_{h_1},H)| + |E_{\widetilde{C}}(P^2_{h_2},H)| - |E_{\widetilde{C}}(b,H)| \leq k+1+k+1 = 2k+2.          \]
By inequality (\ref{eq9}), we have
\begin{equation*}
\widetilde{S} = |E_{\widetilde{C}}(C,H)| \leq l-6-3+3+2k+2+0=l+2k-4.
\end{equation*}
The Claim follows.$\hfill\square$\\

From Claim \ref{claim6}, inequalities (\ref{eq6}) (\ref{eq7}) and (\ref{eq8}), we can deduce that         
\begin{equation*}
\begin{split}
k \delta &\leq \widetilde{S} + S + \widetilde{T} + T +k(k-1) \\
         &\leq l+2k-4 + l+(k-1)(n-l-k)+k(k-1) \\
         &= l+2k-4 + k(n-l)+2l-n.
\end{split}
\end{equation*}
Therefore, we have $k(n-l- \delta +2) \geq n-3l+4$.
Since $l \geq \frac{n+12}{3}$ from Claim \ref{claim1} and $\delta \geq \frac{2n}{3}+1$, we have $n-l- \delta +2 <0$. Then from Claim \ref{claim2}, we have
\[
3(n-l- \delta +2) \geq k(n-l- \delta +2) \geq n-3l+4,
\]
which implies that  $\delta \leq \frac{2n+2}{3}$, a contradiction. We complete the proof of Theorem \ref{thmme}.  $\hfill\square$\\

\section*{Acknowledgment}
This research was supported by National Natural Science Foundation of China under grant numbers 11971445 and 12171440.

\nocite{*}
\bibliographystyle{abbrvnat}

\bibliography{mybib}

\end{document}